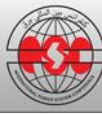

# Optimal Two-Stage Programming for Integration of PHEV Parking Lots in Industrial Microgrids

Farhad Samadi Gazijahani[1*], Javad Salehi[2]

[1,2]Azarbaijan Shahid Madani University
Department of Electrical Engineering
Tabriz, Iran
[1]f.samadi@azaruniv.ac.ir, [2]j.salehi@azaruniv.ac.ir

*Abstract*—With the progressive exhaustion of fossil energy and the enhanced awareness of environmental protection, more attention is being paid to plug in hybrid electric vehicles (PHEV). Inappropriate siting and sizing of PHEV parking lots (PL) could have negative effects on the development of PHEV, the layout of the city traffic network, and the convenience of PHEVs' drivers as well as lead to an increasing in the network losses and a degradation in voltage profiles at some nodes. Given this background, this paper aims to allocate PLs in Industrial Microgrids (IMG) with the objective of minimizing system costs including investment cost, power loss and scheduling cost as possible objectives. A two-stage model has been designed for this purpose. The optimal siting and sizing of PLs in order to minimize the investment cost of PLs is performed in the first stage. At the second stage, the optimal PHEV scheduling problem is solved considering market interactions to provide profit to the PL owner with taken into account various network constraints. Conclusions are duly drawn with a realistic example.

I. INTRODUCTION

*A. Motivation and Aims*

Integrating large numbers of plug in electric vehicles (EVs) into the power grid while simultaneously reducing their impacts and those of uncontrollable renewable energy sources is a major goal of vehicle-to-grid (V2G) systems [1]. V2G is defined as the provision of energy and ancillary services, such as regulation or spinning reserves, from an EV to the grid. This can be accomplished by discharging energy through bidirectional power flow, or through charge rate modulation with unidirectional power flow [2]-[4]. Through V2G, EV owners can produce revenue while their cars are parked which can provide valuable economic incentives for EV ownership. Utilities can also benefit significantly from V2G by having increased system flexibility as well as energy storage for intermittent renewable energy sources such as wind. In order to participate in energy markets, the V2G capabilities of many EVs are combined by aggregators and then bid into the appropriate markets [5]-[8]. An aggregator may be the utility into which the EVs are plugged or a third-party business.

Adequate and proper incentive mechanisms for massive use of EVs, e.g. in the form of monetary incentives to take part in a V2G mode, will lead to sustainable developments in the future grids [9]. Nevertheless, an increasing penetration level of EVs in such systems will materialize if accompanied by sufficient progress for example in the regulatory aspects, and certain requirements are fulfilled [10]. These requirements include communication infrastructures, control and management schemes as well as agents to coordinate their efficient operations [11]. Otherwise, the integration of EVs could lead to undesirable consequences such as technical issues in the systems. Fig. 1 depicts the role of EVs in relation to the new emerging smart MGs paradigm, and illustrates the benefits of EV technology in this context.

*B. Literature Review*

With the continuous technological advances, there are many concepts associated with EVs [12]. The grid-to-vehicle (G2V) and V2G operation modes allow EVs to interact with the grid, either by selling or purchasing power in different periods [13]-[16]. This subject has been widely researched recently. For example, the participation of EVs in the V2G operation mode is reported in [17]-[20]. In [21]-[23], the participation of EVs in the regulation up and down markets are considered. Authors in [24] address EVs' participation in the spinning reserve markets. In addition, the real power systems of Germany and Singapore, which integrate EVs, are described in [25]. Furthermore, the strategic behavior of electric vehicles in their charging models is addressed in [26]-[28]. Although in [26] a reliability cost evaluation model is proposed for a distribution system with both wind generation and PHEVs, it does not consider the V2G mode of PEVs in the reliability study. Reference [27] has studied the real-time coordinated operation of PHEVs in order to minimize distribution network effects, including voltage and loss. However, the study in [28] is more concentrated on the load management aspect of the coordinated charging of PHEV. As a result, this approach is mostly used by the PHEV aggregator. In [29], both V2G and G2V modes of PHEV are studied at different penetration levels to reach acceptable bus voltages and power loss in the grid.

*C. Contributions*

Although the planning of EVs charging stations has been studied in the literature, to the best of our knowledge, the mutual impacts of the planning and operation objectives of EV parking lots on the IMGs have not yet been addressed. Hence, the contributions of this paper can be summarized as follows:

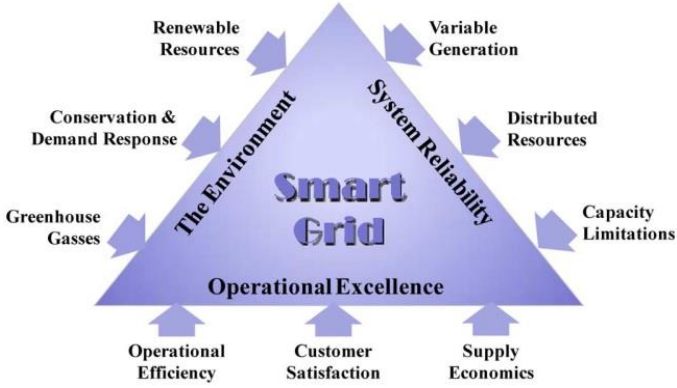

Fig. 1. The role of PHEV in smart MGs [30].

- Developing an innovative model to optimize both planning and operation costs of PHEV parking lots with considering technical constraints of IMGs.
- Proposing an integrated two-stage framework that determines the optimal location and size of PLs at the first stage. Then, this planning is subject to optimal scheduling objectives in order to maximize the profit of PLs owners at the second stage.
- Considering the integrated behavior of PLs through PHEVs' arrivals and departures and also PLs interaction with energy and reserve markets.
- Applying operational planning of PLs with considering various distributed generations.

*D. Paper Organization*

The rest of the paper is organized as follows. Section II presents the two-stage problem to formulate the model of an PHEV parking lot. Section III presents the IMG model featuring RESs. IMG modelling is illustrated in section IV and Section V is devoted to the numerical results of a case study, and some concluding remarks are given in Section VI.

## II. MATHEMATICAL MODEL

This paper proposed a two-stage model to allocate PLs in IMG with the objective of minimizing system costs including investment cost, power loss and scheduling cost as possible objectives. The optimal siting and sizing of PLs in order to minimize the investment cost of PLs is performed in the first stage. At the second stage, the optimal PHEV scheduling problem is solved considering market interactions to provide profit to the PL owner with taken into account various network constraints.

*A. Driving Patterns for EV*

Driving and plug-in patterns of private cars are fairly predictable as many drivers have a fixed schedule of working hours and leisure time activities. This is particularly true for the EVs that are parked for most hours of the day. Statistical data on transportation behavior of such EVs are used to obtain driving patterns [30]. Each pattern represents the distance driven for each hour of the daily driving [31]. In this paper, it is assumed that the categorizes the EVs is clustered into 6 fleets of similar 24-hour driving patterns.

*B. First stage*

In this stage, the objective function is defined as the minimization of the total investment costs associated with PHEV parking lots to be planned. The mathematical model of first stage can be formulated as below:

$$\min_{stage\,1} f^{Up} = \sum_{t=1}^{T} \frac{d(1+d)^t}{(1+d)^t - 1} \left[ \sum_{i=1}^{N_{PL}} x_i^{PL} C_{PL.i}^{IV} S_{PL.i}^R + C_{PL.i}^{OM} S_{PL.i}^R \right.$$

$$\left. + \sum_{l=1}^{N_l} FOR_l P_{l.t}^D C_t^{IL} \right] \quad (1)$$

$$S_{PL.i} \leq S_{PL.i}^{max} \qquad \forall i \in N_{PL} \quad (2)$$

$$I_l \leq I_l^{max} \qquad \forall l \in N_{line} \quad (3)$$

$$V_n^{min} \leq V_n \leq V_n^{max} \qquad \forall n \in N_{bus} \quad (4)$$

$$P_{PL.i}^{min} \leq \sum_{i=1}^{N_{PL}} P_{PL.i} \leq P_{PL.i}^{max} \qquad \forall i \in N_{PL} \quad (5)$$

Where $N_{PL}$ is the number of PHEV parking lots in the IMG concerned, $C_{PL}^{IV}$ and $C_{PL}^{OM}$ are, respectively, the investment cost, as well as operation and maintenance cost of the $i$th PL; $S_{PL.i}^R$ denotes the installed rated size of $i$th PL and also $d$ is the discount rate and used to transform the future cost to the present value; and $T$ is the number of years included in the planning horizon time, which is considered 4 years. Also, $FOR_l$ is forced outage rate of line $l$ and $C_t^{IL}$ is the interrupted cost. The objective function of first stage includes three parts. The first term is the investment cost of PLs; the second terms is the operation and maintenance (O&M) cost of PLs and finally the third part shows the Energy Not Supplied (ENS) cost. The constraints considered here include the equality and inequality ones. Constraint (2) shows the permitted maximal transformer capacity limit of the $i$th PL. Constraints (3) and (4) guarantees that the current $I_l$ and voltage $V_n$ of system must be located in the acceptable regions. The permitted minimal and maximal charging power limits of PHEV charging stations are illustrated in constraint (5).

*C. Second stage*

At the second stage, the optimal scheduling of PLs in an IMG is performed in order to maximize the profit of PLs operator's viewpoint considering various network-constrained objectives. In this study, as shown in (6), the profit is obtained through energy and reserve market interactions individual contracts with PHEV owners that use the PLs in a V2G state as well as decreasing in the power losses of system.

$$\max_{stage\,2} f^{Dn} = \sum_{t=1}^{4} IR \sum_{d=1}^{365} \sum_{h=1}^{24} \left\{ \left( \sum_{i=1}^{N_{PL}} P_{PL.i.t.d.h}^S \pi_{PL.i.t.d.h}^S \right. \right.$$

$$\left. \left. - P_{PL.i.t.d.h}^B \pi_{PL.i.t.d.h}^B \right) - \sum_{l=1}^{N_l} R_l I_{l.t.d.h}^2 \right\} \quad (6)$$

$$P^S_{PL.i.t.d.h} \leq P^{max}_{PL.i} \quad (7)$$

$$P^B_{PL.i.t.d.h} \leq P^{max}_{PL.i} \quad (8)$$

$$SOC^{min}_i \leq SOC_i \leq SOC^{max}_i \quad (9)$$

$$IR = \frac{d(1+d)^t}{(1+d)^t - 1} \quad (10)$$

$$SOC_{i.t.d.h} = SOC^{INI}_{i.t.d.h} + (\eta_c P^{CH}_{i.t.d.h} - \eta_d P^{DCH}_{i.t.d.h}) \quad (11)$$

$$P^{CH}_{PL.i} + P^D_{t.d.h} + P^S_{t.d.h} = P^{DCH}_{PL.i} + P^B_{t.d.h} \quad (12)$$

The aim of second stage is to maximize the profit of PL's owners during planning horizon time. The equation (6) depicts the profit based objective function, where includes three terms. The two first terms show the sold and purchased power of PLs and third term is the power losses. Constraints (7) and (8) demonstrate the maximum allowed sold and purchased power through PLs, respectively. Equations (9) to (11) show the maximum and minimum state of charge (SOC) of PHEVs, interest rate and the SOC of PLs in each operation interval, respectively. Furthermore, the constraint (12) shows the power mismatch of IMG which determines that at each time, sum of total generated powers must be equal to total consumed powers. $P^S_{PL.i.t.d.h}$ and $\pi^S_{PL.i.t.d.h}$ denote the sold power and sold price by $i$th PL at $t$th year, $d$th day and $h$th hour, respectively. Also, $P^B_{PL.i.t.d.h}$ and $\pi^B_{PL.i.t.d.h}$ show the purchased power and purchased price by $i$th PL at $t$th year, $d$th day and $h$th hour, respectively. $R_l$ is the resistance of line $l$ which can be obtained from data of utilized network and $I_{l.t.d.h}$ is the current of line $l$ at year $t$, day $d$ and hour $h$. Furthermore, $\eta_c$ and $\eta_d$ are charged and discharged efficiencies of PHEVs, respectively. Also, $P^{CH}_{i.t.d.h}$ and $P^{DCH}_{i.t.d.h}$ denote the charged and discharged power of $i$th PHEV at $t$th year, $d$th day and $h$th hour.

Note that in the PLs, the SOC is a measure of the amount of energy stored in the PHEVs. It is similar to the fuel gauge in conventional internal combustion cars. In this paper, SOC refers to the percentage of energy remained in the battery when PHEV arrives home, after daily trips. So, PHEV can operate in blended mode in which the internal combustion engine helps the electric motor provide the required energy to run the vehicle. Therefore, in charge depleting mode, either the entire or a fraction of the required energy is supplied by the battery. To cover all the possible PHEV operations, a factor for each PHEV is defined. This factor represents the percentage of distance that PHEV drives in the electric mode. Accordingly, the SOC of a PHEV would be as below [32].

$$SOC = \begin{cases} \left(1 - \dfrac{\psi \times \phi}{AER}\right) & \psi \times \phi \leq AER \\ 0 & \psi \times \phi \geq AER \end{cases} \quad (13)$$

Where $\phi$ is the total driven distance, $(\psi \times \phi)$ is the distance driven in the electric model. The battery will be empty if reaches all electric range (AER). The amount of energy required to charge the PHEV battery can be calculated as follows [33].

$$E_C = \left(1 - \frac{SOC}{100}\right) \times C \quad (14)$$

$$C = \beta \times AER \quad (15)$$

$$E_g = \frac{E_C}{\xi} \quad (16)$$

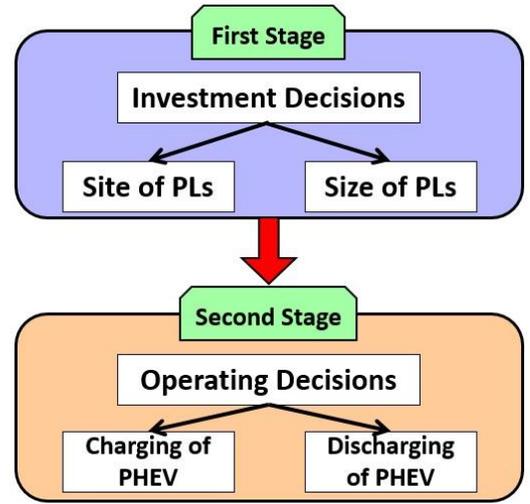

Fig. 2. Flowchart of proposed two-stage programming for PLs planning.

where $\xi$ is the efficiency of charging PHEV battery and $E_g$ is the actual energy which should be transferred from the grid to charge the battery. The efficiency of PHEV is considered to be 90%. Also, $C$ denotes the usable capacity of PHEV battery and $\beta$ is the electrical energy consumption per mile [37].

### III. SOLUTION METHODOLOGY

In this paper, the proposed two-stage model for PLs planning problem has been formulated as Mixed Integer Linear Programming (MILP) for both stages and subsequently solved using Linear Programming (LP) in GAMS software by CPLEX solver. Linear programming is considered a revolutionary development that permits us to make optimal decisions in complex situations. Linear programming deals with the problem of optimizing (minimizing or maximizing) a linear function of $n$ variables subject to equality and/or inequality linear constraints. In other words, a linear programming problem consists of finding the optimum of a linear function in a set that can be written as the intersection of a finite number of hyperplanes and half spaces in $\mathcal{R}^n$. Although several other methods have been developed over the years for solving LP problems, the simplex method continues to be the most efficient and popular method for solving general LP problems. The general LP problem for first and second stages can be stated in the following standard forms:

$$Min\ f^{Up}\ Z(X) = C^T X \quad (17)$$
Subject to: (2) to (5) $\quad (18)$
$$Max\ f^{Dn}\ Z(Y) = D^T Y \quad (19)$$
Subject to: (7) to (16) $\quad (20)$

Fig. 2 shows the proposed two-stage model for optimal planning of PL. With respect to this figure, it can be seen that by applying proposed method, the planning and operation objectives of PLs integration problem could be handled coordinately.

Table I. Technical data of elements of IMG.

| Elements | Number | Min (kW) | Max (kW) |
|---|---|---|---|
| CHP | 2 | 200 | 500 |
| WT | 2 | 100 | 300 |
| Switch | 6 | 1, 5, 10, 14, 26, 31 | |

## IV. IMG M ODELLING

The test system is the 37-bus system [34] IMG consisting of 2 factories with Combined Heat and Power (CHP) systems coupled with Wind Turbine (WT). The Distributed Energy Resources (DER) data is given in Table I [35]. All factories cooperate in generating electricity. However, only neighboring factories are allowed to participate in acquiring the required heat. IMGs rely on CHP systems to facilitate energy-efficient power generation by capturing the waste heat. These systems maintain the heat acquired from power generation and utilize it for domestic and industrial heating purposes [36]. IMGs can be connected or disconnected from the upstream network. In the stand-alone mode, IMGs must generate their own required energy to feed the electric loads through the cooperation of all DG units. In the grid-connected mode, IMGs are permitted to purchase some of their electric needs from the upstream network or even sell electricity to upstream network in some hours of the day. On the other hand, due to the existing distances between factories, only the ones in the vicinity of each other can cooperate to procure the thermal needs. Of course, some factories may not have thermal requirements. The daily load curve as shown in Fig. 3 is used for operation scheduling of PHEVs during planning horizon time. It should be mentioned that in this paper, the uncertainty of WT generation is taken into account through possibility approach which considers the expected value of uncertain parameters obtained from Weibull distribution [32]. The expected value is most likely to occur. The Weibull distribution of WTs is demonstrated in Fig. 8.

## V. S IMULATION AND R ESULTS

The IMG is an office block and the involved PHEVs are private cars and official vehicles. All the PHEVs have the same batteries. Some assumptions about the PHEV system are given as follows. The PHEV charging and discharging time meets the owners' habits [37]. The charging and discharging time of the official EVs (OEVs) is assumed to start after work at 16:00 to 7:00 (next day). In order to meet the needs of the next day, OEVs need to get charged every day and the disconnection SOC is not less than 98% at the end time. The charging and discharging of the private EVs (PHEVs) can take place during the working hours or after work at home. Thus, the charging and discharging time of the PEVs in the office block is assumed at 8:00–15:00 and the disconnection SOC is not less than 45% to ensure the normal running after work.

The energy consumption parameter of daily distance driven of the PHEVs is considered to follow normal PDF. In this paper, it is assumed that the energy consumption parameter of daily one-way distance driven for PEVs is between 0 and 0.25; it is assumed that the energy consumption parameter of OEVs is between 0 and 0.4. The MG can provide a backup battery to the owner and accept the battery as the echelon-use battery when the EV decline in cell performance [38].

In this paper, two case studies are investigated to evaluate the proposed model. On the first approach, the Pay as Bid pricing model is investigated to examine the individual interaction of PHEV with the aggregator [39]. In the second approach, the cross effect of the resources in their market participation is investigated through uniform pricing. The problem is modeled as an MILP problem and implemented in GAMS using CPLEX solver [40].

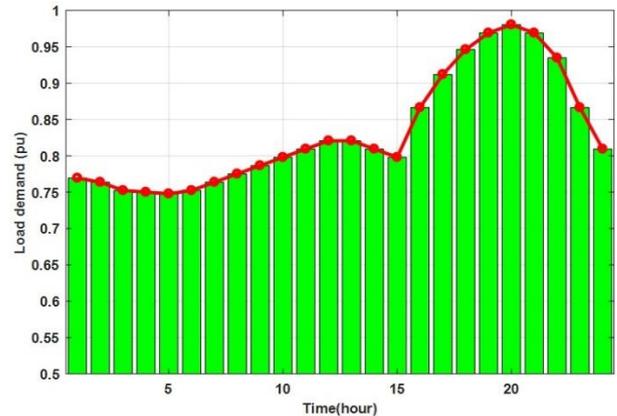

Fig. 3. Forecasted IMG daily electrical load.

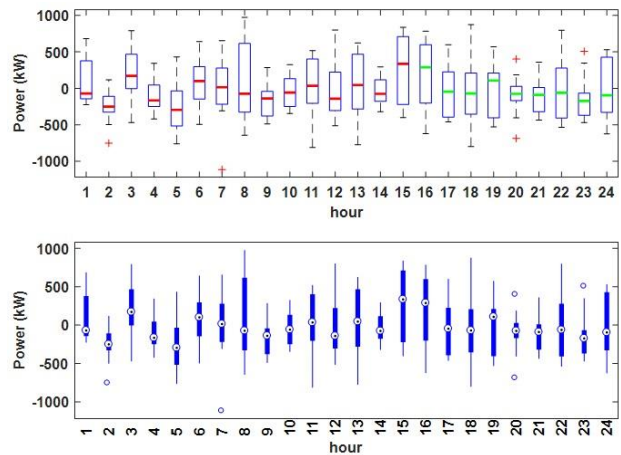

Fig. 4. Box plot of installed PLs in one day.

Fig. 4 illustrates the charging and discharging of PLs in the daily scheduling. As can be observed in Fig. 4, the availability of more energy resources in the grid lets the parking lot purchase a higher quantity of energy from the energy market. This is more evident in case of without PL integration, which consider wind power generators in the distribution networks. Consequently, the amount of energy sold back to the grid is also higher. This reveals that the contribution of the parking lot in the energy market in a V2G mode in cases without DER units is higher than that of case of DER units. Table II shows the obtained capacity for PLs in the IMGs. According to this table, it can be seen that six PLs integrated into IMG. In this paper, all the stations in the PLs are the same and are quick charging stations with a charging rate of 11 kW per hour [41].

Table II. Site and Size of installed PLs on the IMG.

| PLs | Site | Size [kW] |
|---|---|---|
| PL-1 | 2 | 250 |
| PL-2 | 17 | 200 |
| PL-3 | 25 | 200 |
| PL-4 | 22 | 225 |
| PL-5 | 11 | 275 |
| PL-6 | 34 | 250 |

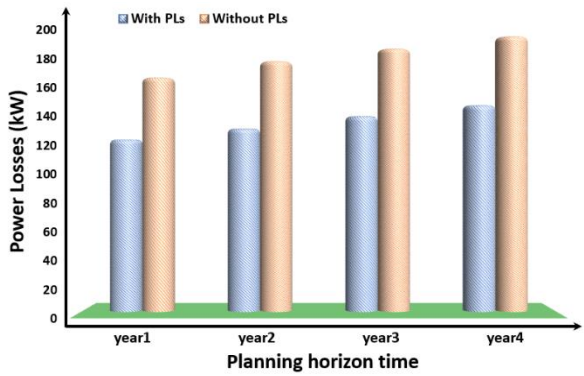
Fig. 5. Power losses of IMG with and without PLs.

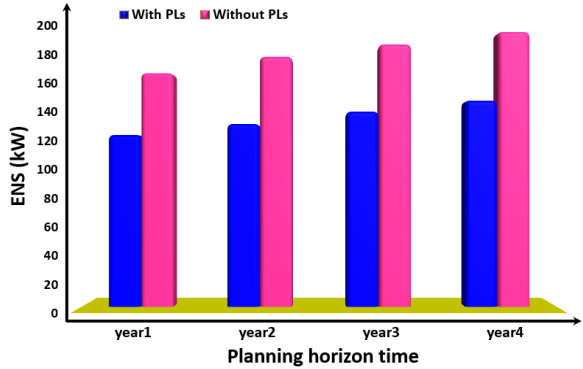
Fig. 6. ENS value of IMG with and without PLs.

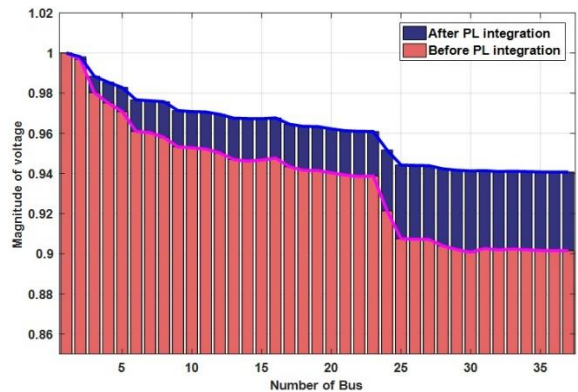
Fig. 7. Voltage profile of IMG before and after integration of PLs.

the PHEVs departure from the PL increases. As a result, in order to meet the PHEV's preferences, the charging of PL is limited. For the reserve provision, except where the reserve price faces a spike at hour 15, in other hours the price is almost equal to the marginal price of PL for providing reserve.

Table III. Comparison between with and without PL cases.

| Costs | Without PLs | With PLs |
|---|---|---|
| Voltage deviations | 9.64% | 4.14% |
| Power losses | 745 kW | 486 kW |
| Total ENS value | 528 kWh | 346 kWh |

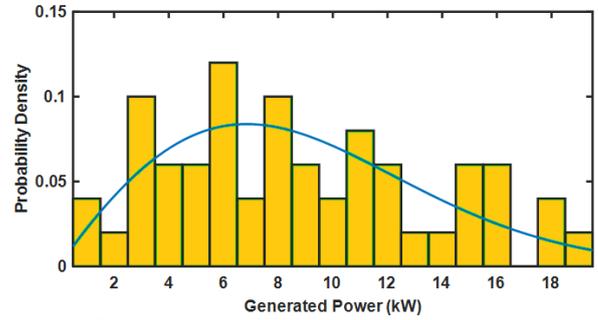
Fig. 8. The Weibull distribution of WT generation.

Fig. 5 and Fig. 6 illustrate the power losses and ENS of IMG in both cases. Regarding to this figure, it is clear that by utilizing PLs in the IMG the power losses of system as well as ENS cost will be decreased significantly. Furthermore, the voltage profile of proposed network is demonstrated in Fig. 7. It is obvious that PLs can appropriately be improved the voltage deviations. The voltage profile of load points at peak times, which parking lots deliver power to the distribution system is shown in this figure. As can be observed from Fig. 7, there are improvements in voltage profile of some buses in the presence of V2G power. Table III compares the technical characteristics of system in cases of with and without PLs integration. As regards to this table, it can be concluded that PLs could be improved various aspects of IMG such as ENS value. This table compares the objective function values under the two cases of IMG energy scheduling approaches. As most vehicles are expected to park for longer than the actual charging time, the constrained charging scheme can reduce the operational cost of the IMG by taking full advantage of timing flexibility.

In this paper, the PL as the main concern of the study changes its behavior based on its trade with the PEV owners and the aggregator. As a result, the tariffs that are implemented to the PHEVs can significantly change the strategy of the PL in the market. The variation of the behavior also leads to different levels of profit gain for the PL and aggregator [43]. In this study, PL is a complicated resource in the system which can act as a flexible demand and as a resource as well. Therefore, the aggregator can benefit the most from the PL's potential to act as the flexible load. However, the aggregator needs to manage the market wisely to encourage the PL to show more flexibility. The proposed two-stage model has been implemented on the modified IEEE 37-bus distribution system as shown in Fig. 9 which used to verify the proposed model. The microgrid DGs, such as CHP and WT units as well as installed PLs, which are obtained from proposed optimization process are shown in this figure.

However, both the PL's behavior and the price are propelled to the equilibrium price as in this price the optimum profit is obtained. During the early hours of the day (hours 1-9) the PL starts to charge the PHEVs in the PL because the energy price is low [42]. The PL can make profit from selling energy to the PHEVs, however the preferences of PHEVs on requiring a fixed amount of departure SOC limits the charging behavior of the PL. Meanwhile, the aggregator wants to increase its profit from selling energy to the PL; as a result, it will encourage the PL to charge its PHEVs by increasing the price of reserve at hours 10 and 15 (see Fig. 4). The price of reserve is increased by the aggregator so that the PL will be motivated for charging; however, the preferences of the PHEVs limit the maximum charging of PL. In fact, noting Fig. 4, it is shown that the PHEVs are charged almost the same as their minimum requirement of departure SOC. The reason is that from hour 15,

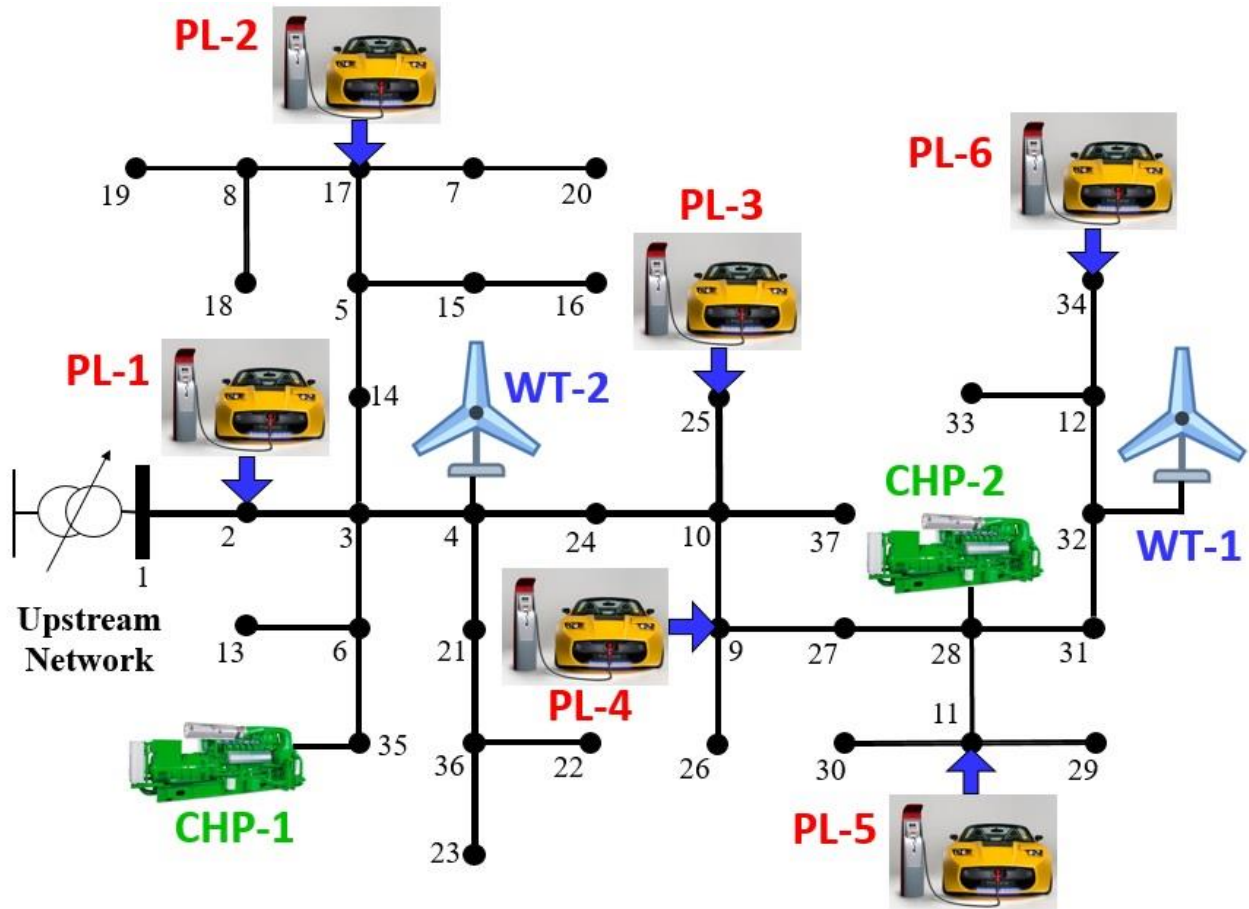

Fig. 9. Optimal allocation of PLs in the applied IMG.

## VI. Conclusion

To solve the optimal planning problem of EV parking lots, a method combining the two-stage screening method which has been formulated as MILP model. The two-stage screening method with the economic and technical factors and the service radius of EV charging stations considered is first presented to identify the optimal sites and sizes for EV charging stations. Then, a mathematical model for the optimal scheduling of EV charging stations is developed, and solved by CPLEX solver in GAMS. Finally, simulation results of the IEEE 37-node test feeder demonstrate that the developed model and method cannot only attain the reasonable planning scheme of EV charging stations, but also reduce the network loss and improve the voltage profile. The results show that a PL, due to its nature as a charging station, will behave more likely like a load in the system. However, in certain situations, the V2G mode can be used and the PL will act as a resource in the system. As a result, by optimal arbitrage of PL in power market can provide more benefit for operator as well as causes to improve the technical specifications of system such as power loss and ENS. Also, more revenue from PHEV owners can be obtained due to a higher SOC that will remain in the PHEV batteries. Regarding network-constrained objectives, despite the low costs of V2G for PL owners, microgrid operator can profit significantly from the presence of PLs in the system.